\documentclass[reqno,11pt]{amsart}
\usepackage{amscd,amssymb,verbatim}

\numberwithin{equation}{section}
\theoremstyle{plain}

\renewcommand{\widetilde}{\tilde}

\newcommand\alp{\alpha}         

\newcommand\gam{\gamma}         \newcommand\Gam{\Gamma}
         
\newcommand\eps{\varepsilon}

\newcommand\iot{\iota}

\newcommand\lam{\lambda}                \newcommand\Lam{\Lambda}
\newcommand\sig{\sigma}         

\newcommand\ome{\omega}         

\newcommand\calA{{\mathcal{A}}}

\newcommand\calH{{\mathcal{H}}}

            \newcommand\bfB{{\mathbf B}}

\newcommand\RR{\mathbb{R}}

\newcommand\CC{\mathbb{C}}

\newcommand\sdp{\times \hskip -0.3em {\raise 0.3ex
\hbox{$\scriptscriptstyle |$}}} 

\newcommand\Dom{\operatorname{Dom}}

\newcommand\ind{\operatorname{ind}}

\newcommand\Ker{\operatorname{Ker}}

\newcommand\rk{\operatorname{rk}}

\newcommand\oj{{\overline{j}}}

\newcommand\tilB{{\widetilde{B}}}

\newcommand\tilF{{\widetilde{F}}}

\newcommand\tilW{{\widetilde{W}}}

\renewcommand{\>}{\rangle}
\newcommand{\<}{\langle}

\theoremstyle{plain}
\newtheorem{Thm}[subsection]{Theorem}
\newtheorem{Cor}[subsection]{Corollary}
\newtheorem{Lem}[subsection]{Lemma}
\newtheorem{Prop}[subsection]{Proposition}
\newtheorem{Conjec}[subsection]{Conjecture}

\theoremstyle{definition}
\newtheorem{Def}[subsection]{Definition}

\theoremstyle{remark}

\newtheorem{Rem}[subsection]{Remark}

\newif\ifShowLabels
\ShowLabelstrue
\newdimen\theight
\def\TeXref#1{%
        \leavevmode\vadjust{\setbox0=\hbox{{\tt
                \quad\quad  {\small \textrm #1}}}%
        \theight=\ht0
        \advance\theight by \lineskip
        \kern -\theight \vbox to
        \theight{\rightline{\rlap{\box0}}%
        \vss}%
        }}%

\newcommand{\refs}[1]{Section ~\ref{S:#1}}
\newcommand{\refss}[1]{Subsection ~\ref{SS:#1}}
\newcommand{\reft}[1]{Theorem ~\ref{T:#1}}
\newcommand{\refl}[1]{Lemma ~\ref{L:#1}}
\newcommand{\refp}[1]{Proposition ~\ref{P:#1}}
\newcommand{\refc}[1]{Corollary ~\ref{C:#1}}

\newcommand{\refe}[1]{\eqref{E:#1}}

\newenvironment{thm}[1]%
        { \begin{Thm} \label{T:#1}  \ifShowLabels \TeXref{T:#1} \fi }%
        { \end{Thm} }

\renewcommand{\th}[1]{\begin{thm}{#1} \sl }
\renewcommand{\eth}{\end{thm} }

\newenvironment{lemma}[1]%
        { \begin{Lem} \label{L:#1}  \ifShowLabels \TeXref{L:#1} \fi }%
        { \end{Lem} }
\newcommand{\lem}[1]{\begin{lemma}{#1} \sl}
\newcommand{\elem}{\end{lemma}}

\newenvironment{propos}[1]%
        { \begin{Prop} \label{P:#1}  \ifShowLabels \TeXref{P:#1} \fi }%
        { \end{Prop} }
\newcommand{\prop}[1]{\begin{propos}{#1}\sl }
\newcommand{\eprop}{\end{propos}}

\newenvironment{corol}[1]%
        { \begin{Cor} \label{C:#1}  \ifShowLabels \TeXref{C:#1} \fi }%
        { \end{Cor} }
\newcommand{\cor}[1]{\begin{corol}{#1} \sl }
\newcommand{\ecor}{\end{corol}}

\newenvironment{conjec}[1]%
        { \begin{Conjec} \label{Conj:#1}  \ifShowLabels \TeXref{C:#1} \fi }%
        { \end{Conjec} }
\newcommand{\conj}[1]{\begin{conjec}{#1} \sl }
\newcommand{\econj}{\end{conjec}}

\newenvironment{defeni}[1]%
        { \begin{Def} \label{D:#1}  \ifShowLabels \TeXref{D:#1} \fi }%
        { \end{Def} }
\newcommand{\defe}[1]{\begin{defeni}{#1} \sl }
\newcommand{\edefe}{\end{defeni}}

\newenvironment{remark}[1]%
        { \begin{Rem} \label{R:#1}  \ifShowLabels \TeXref{R:#1} \fi }%
        { \end{Rem} }
\newcommand{\rem}[1]{\begin{remark}{#1}}
\newcommand{\erem}{\end{remark}}

\newcommand{\eq}[1]%
        { \ifShowLabels\newline \TeXref{E:#1} \fi
           \begin{equation} \label{E:#1} }
\newcommand{\eeq}{\end{equation}}

\newcommand{\prf}{ \begin{proof} }
\newcommand{\eprf}{ \end{proof} }
\newcommand{\Label}[1]{\label{#1}  \ifShowLabels \TeXref{#1} \fi }

\ShowLabelsfalse

\renewcommand{\d}{\text{\( \partial\)}}

\renewcommand{\b}{\bullet}

\newcommand{\A}{\calA}

\renewcommand{\i}{\sqrt{-1}\, }
\newcommand{\mB}{\bfB^{\text{mod}}}

\renewcommand{\H}{\bfB_a^2}

\begin{document}
\title{New proof of the cobordism invariance of the index}
\author{Maxim Braverman}
\address{Department of Mathematics\\
        Northeastern University   \\
        Boston, MA 02115 \\
        USA
         }
\email{maxim@neu.edu}
\begin{abstract}
We give a simple proof of the cobordism invariance of the index of
an elliptic operator. The proof is based on a study of a
Witten-type deformation of an extension of the operator to a
complete Riemannian manifold. One of the advantages of our
approach is that it allows to treat directly general elliptic
operator which are not of Dirac type.
\end{abstract}
\maketitle \vspace{-0.8cm}
\section{Introduction}\Label{S:setting}

Recently several simple proofs of the cobordism invariance of the
index were established, cf. \cite{Higson91},
\cite[Th.~6.2]{Lesch93}, \cite{Nicol97}. In this note we present
still another proof of this fact. Unlike other authors we don't
impose any restrictions on the dimension of the manifold and don't
assume that our operator is of Dirac type.

\subsection{The setting}\Label{SS:setting}
Let $E^+,E^-$ be Hermitian vector bundles over a closed Riemannian
manifold $M$. Let $A^+:C^\infty(M,E^+)\to C^\infty(M,E^-)$ be an
elliptic differential operator. Let $A^-:C^\infty(M,E^-)\to
C^\infty(M,E^+)$ be the formal adjoint of $A^+$ and consider the
operator
\[ A \ := \
 \begin{bmatrix}
  0&A^-\\A^+&0
 \end{bmatrix}
 : \ C^\infty(M,E^+\oplus E^-) \ \to \ C^\infty(M,E^+\oplus E^-).
\]
This operator is essentially self-adjoint and we denote by the
same letter $A$ its extension to a self-adjoint operator acting on
the space $L^2(M,E^+\oplus{E^-})$ of square-integrable sections.

Suppose now that $M$ is a boundary of a Riemannian manifold $W$,
which is isometric near the boundary to the cylinder
$U=M\times(-\eps,0]$. Let $F$ be a Hermitian vector bundle over
$W$, whose restriction to $U$ is isomorphic to the lift of
$E^+\oplus{E^-}$.

\th{cobinv}
Assume that there exists a self-adjoint hypo-elliptic differential
operator $B:C^\infty(W,F)\to C^\infty(W,F)$, which near the
boundary takes the form
\[
    B \ = \ \gam\frac{\d}{\d t} \ + \ A,
\]
where $t$ is the normal coordinate and $\gam$ is a skew-adjoint
bundle map independent of $t$ such that $\gam|_{E^\pm}=\pm\i$.
Then the index \/
 \(\displaystyle
 \ind{A} \ := \ \dim\Ker{A^+} \ - \ \dim\Ker{A^-} \ = \ 0.
 \)
\eth

\subsection{The plan of the proof}\Label{SS:plan} Let $\tilW$ denote
the complete non-compact Riemannian manifold obtained from $W$ by
attaching the semi-infinite cylinder $M\times[0,\infty)$  to the
boundary. We extend the bundle $F$ and the operator $B$ to $\tilW$
in the obvious way.

Consider the exterior algebra
$\Lam^\b\CC=\Lam^0\CC\oplus\Lam^1\CC$. It has two (anti)-commuting
actions $c_L$ and $c_R$ (left and right action) of the Clifford
algebra of $\RR$, cf. \refss{wE}. Set $\tilF=F\otimes\Lam^\b\CC$
and consider the operator
\eq{tilA}
    \tilB \ := \ \i B\otimes c_L(1): \, C^\infty(\tilW,\tilF)
     \ \to \ C^\infty(\tilW,\tilF).
\end{equation}

Let $p:\tilW\to\RR$ be a map, whose restriction to
$M\times(1,\infty)$ is the projection on the second factor, and
such that $p(W)=0$ (see \refss{wE} for a convenient choice of this
function). For any $a\in\RR$, consider the operator $\bfB_a:=\tilB
- 1\otimes{}c_R((p(t)-a))$. Then (cf. \refl{bfB2})
\eq{bfA2}
    \bfB_a^2 \ = \ B^2\otimes1 - R + |p(x)-a|^2,
\end{equation}
where $R:\Gam(\tilW,\tilF)\to\Gam(\tilW,\tilF)$ is a bounded
operator.

Set $\ind\bfB_a:=\dim\Ker\bfB_a^+-\dim\Ker\bfB_a^-$, where
$\bfB_a^\pm$ denote the restriction of $\bfB_a$ to the spaces
$F\otimes\Lam^0\CC$ and $F\otimes\Lam^1\CC$ respectively. It
follows  from \refe{bfA2} that $\ind\bfB_a=0$ for $a\ll 0$ and, if
$a\gg0$, then all the sections in $\Ker\bfB_a^2$ are concentrated
on the cylinder $M\times(0,\infty)$, not far from $M\times\{a\}$
(this part of the proof essentially repeats the arguments of
Witten in \cite{Witten82}). Hence, the calculation of
$\Ker\bfB_a^2$ is reduced to a problem on the cylinder
$M\times(0,\infty)$. It is not difficult now to show that
$\ind\bfB_a=\ind{A}$ for $a\gg0$.

\reft{cobinv} follows now from the fact that $\ind\bfB_a$ is
independent of $a$.

\section{Index of the operator $\bfB_a$}\Label{S:Aa}

\subsection{}\Label{SS:wE}
Let us consider two anti-commuting actions (left and right action)
of the Clifford algebra of $\RR$ on the exterior algebra
$\Lam^\b\CC=\Lam^0\CC\oplus\Lam^1\CC$, given by the formulas
\eq{ClonR}
    c_L(t)\, \ome \ = \ t\wedge\ome \ - \ \iot_t\ome; \qquad
    c_R(t)\, \ome \ = \ t\wedge\ome \ + \ \iot_t\ome.
\end{equation}

We will use the notation of \refss{plan}. In particular, $\tilW$
is the manifold obtained from $W$ by attaching a cylinder,
$\tilF=F\otimes\Lam^\b\CC$ and $\tilB$ is the operator defined in
\refe{tilA}.

Let $s:\RR\to[0,\infty)$ be a smooth function such that $s(t)=t$
for $|t|\ge1$, and $s(t)=0$ for $|t|\le1/2$. Consider the map
$p:\tilW\to\RR$ such that $p(y,t)=s(t)$ for $(y,t)\in
M\times(0,\infty)$ and $p(x)=0$ for $x\in W$. Define the operator
\eq{bfB}
    \bfB_a \ := \ \tilB \ - \ 1\otimes{}c_R((p(x)-a)).
\end{equation}
The same proof as in \cite[Th.~1.17]{GromLaw83}, shows that the
operator $\bfB_a$ is essentially self-adjoint with the initial
domain smooth compactly supported sections. We will denote by
$\bfB_a$ also the extension of this operator to a self-adjoint
operator on the space of square-integrable sections.

\lem{bfB2}
Let $\Pi_i:\tilF\to F\otimes\Lam^i\CC$, \/ $(i=0,1)$ be the
projections. Then
\eq{bfB2}
    \bfB_a^2 \ = \  B^2\otimes1  \ - \ R \ + \ |p(x)-a|^2,
\end{equation}
where $R:\tilF\to\tilF$ is a uniformly bounded bundle map, whose
restriction to $M\times(1,\infty)$ is equal to
$\i\gam(\Pi_1-\Pi_0)$, and whose restriction to $W$ vanishes.
\elem
\prf
Note, first, that $p(x)-a\equiv-a$ on $W$. Thus, since $c_R(a)$
anti-commutes with $\tilB$, we have $\bfB_a^2|_W= \tilB^2|_W+a^2=
B^2\otimes1|_W+a^2$. Hence,  \refe{bfB2} holds, when restricted to
$W$.

We now consider the restriction of $\bfB_a^2$ to the cylinder
$M\times(0,\infty)$. Recall that the function $s:\RR\to[0,\infty)$
was defined in \refss{wE}. Clearly,
\[
    \bfB_a|_{M\times(0,\infty)} \ = \
       \i B\otimes c_L(1)
        \ + \ \i\gam\otimes c_L(1)\frac{\d}{\d t}
         \ + \ \big(\,  s(t)-a\, \big)\, 1\otimes c_R(1).
\]
Since the operators $c_L$ and $c_R$ anti-commute, we obtain
\[
    \bfB_a^2|_{M\times(0,\infty)} \ = \
    \ = \
        B^2\otimes1 \ + \ \i s'\gam\otimes c_L(1)c_R(1) \ + \  |t-a|^2.
\]

Since $c_L(1)c_R(1)=\Pi_1-\Pi_0$, it follows, that \refe{bfB2}
holds with $R=s'\i\gam\big(\Pi_1-\Pi_0\big)$.
\eprf
\lem{discrete}
The spectrum of the operator $\bfB_a$ is discrete.
\elem
\prf
It is well known, cf., for example, \cite[Lemma~6.3]{Sh99}, that
the Lemma is equivalent to the following statement: For any
$\eps>0$ there exists a compact set $K\subset \tilW$, such that if
$u$ is a smooth compactly supported section of $\tilF$, then
\eq{<eps}
    \int_{\tilW\backslash K}\, |u|^2\, d\mu \ < \
      \eps \int_\tilW\, \<\bfB_a^2u,u\>\, d\mu.
\end{equation}
Here, $d\mu$ is the Riemannian volume element on $\tilW$, and
$\<\cdot,\cdot\>$ denotes the Hermitian scalar product on the
fibers of $\tilF$.

Set $V(x)=|p(x)-a|^2-R$. To prove \refe{<eps} note that, since $R$
is a bounded, there exists a compact set $K\subset \tilW$, such
that $V>1/\eps$ on ${\tilW\backslash K}$.
Note, also, that the
first summand in \refe{bfB2} is a non-negative operator. Hence, we
have
\[
    \int_{\tilW\backslash K}\, |u|^2\, d\mu \ < \
      \eps\int_{\tilW\backslash K}\, \<Vu,u\>\, d\mu
      \ \le \  \eps\int_{\tilW}\, \<Vu,u\>\, d\mu \ \le \
        \eps\, \int_{\tilW}\, \<\bfB_a^2u,u\>\, d\mu.
\]
\eprf
Set $\tilF^+:= F\otimes\Lam^0\CC, \ \tilF^-:= F\otimes\Lam^1\CC, \
\bfB_a^\pm:= \bfB_a|_{\Gam(\tilW,\tilF^\pm)}$ and define
\eq{indB}
    \ind\bfB_a \ = \ \dim\Ker\bfB_a^+ \ - \ \dim\Ker\bfB_a^-.
\end{equation}
\lem{indepa}
The index $\ind\bfB_a$ is independent of $a$.
\elem
\prf
From \refe{bfB}, we see that $\bfB_b-\bfB_a=1{\otimes}c_R(b-a)$ is
a bounded operator, depending continuously on $b-a\in\RR$. The
lemma follows now from the stability of the index of a Fredholm
operator, cf., for example, \cite[\S{}I.8]{Sh1}.
\eprf

\lem{depofa}
\(\displaystyle\ind(\bfB_a)=0 \) for all $a\in\RR$.
\elem
\prf
By \refl{indepa}, it is enough to prove the proposition for one
particular value of $a$. But it follows from \refl{bfB2} that, if
$a$ is a negative number such that $a^2>\sup_{x\in\tilW}\|R(x)\|$,
then $\bfB_a^2>0$, so that $\Ker\bfB^2_a=0$.
\eprf

To prove \reft{cobinv} it is enough now to show that $\ind\bfB_a=
\ind{A}$. This is done in two steps: first, in \refs{model}, we
construct a ``model" operator $\mB$ on the cylinder
$M\times(-\infty,\infty)$, whose index is equal to $\ind{A}$.
Then, in \refs{pra-a}, we show that $\ind\bfB_a= \ind\mB$.


\section{The model operator}\Label{S:model}

The bundles $E^\pm$ lift to Hermitian vector bundles over the
cylinder $M\times\RR$, which we will denote by the same letters.
Consider the Hermitian vector bundle $\tilF:=
(E^+\oplus{E^-})\otimes\Lam^\b\CC$ and the operator $\mB:
C^\infty(M\times\RR,\tilF)\to C^\infty(M\times\RR,\tilF)$ defined
by
\[
  \mB \ := \ \i B\otimes c_L(1)
        \ + \ \i\gam\otimes c_L(1)\frac{\d}{\d t}
         \ + \ 1\otimes c_R(t),
\]
where $t$ is the coordinate along the axis of the cylinder. We
refer to $\mB$ as the {\em model operator}, cf. \cite{Sh5}. As in
\refs{Aa}, it is essentially self-adjoint and has discrete
spectrum. We define $\ind\mB$ by \refe{indB}.
\lem{model}
The kernel of the model operator $\mB$ is isomorphic (as a graded
space) to $\Ker(A)$. In particular, $\ind\mB=\ind A$.
\elem
\prf
The same calculations as in the proof of \refl{bfB2}, show that
\eq{Dmodpm}
    (\mB)^2|_{\Gam(M\times\RR,E^\pm\otimes\Lam^\b\CC)} \ = \
    A^2\otimes1 \ + \
    1\otimes\, \Big(\, -\frac{\d^2}{\d t^2}\pm(\Pi_1-\Pi_0)+t^2
    \, \Big).
\end{equation}
Both summands in the right hand side of \refe{Dmodpm} are
non-negative. Hence, the kernel of $(\mB)^2$ is given by the
tensor product of the kernels of these operators.

The space $\Ker\big(-\frac{\d^2}{\d t^2}+\Pi_1-\Pi_0+t^2\big)$ is
one dimensional and is spanned by the function
$\alp^+(t):=e^{-t^2/2}\in \Lam^0\RR$. Similarly,
$\Ker\big(-\frac{\d^2}{\d t^2}+\Pi_0-\Pi_1+t^2\big)$ is one
dimensional and is spanned by the one-form
$\alp^-(t):=e^{-t^2/2}ds$, where we denote by $ds$ the generator
of $\Lam^1\CC$. It follows that

 \(\displaystyle\hskip1.3in
    \Ker(\mB)^2|_{\Gam(M\times\RR,E^\pm\otimes\Lam^\b\CC)} \ \simeq \
    \Big\{\, \sig\otimes\alp^\pm(t): \
    \sig\in\Ker A^2|_{\Gam(M,E^\pm)}\, \Big\}.
 \)
\eprf

\subsection{}\Label{SS:DaV}
Let $T_a:M\times\RR\to M\times\RR, \ T_a(x,t)=(x,t+a)$ be the
translation and consider the pull-back map
$T_a^*:\Gam(M\times\RR,\tilF)\to \Gam(M\times\RR,\tilF)$. Set
\[
    \mB_a \ := \ T_{-a}^*\circ\mB\circ T^*_a \ = \
        B\otimes1 \ - \ 1\otimes c_R\big(\, t-a\, \big)
\]
Then \/ $\ind\mB_a=\ind\mB$, for any $a\in\RR$.

\section{Proof of \reft{cobinv}}\Label{S:pra-a}

If $A$ is a self-adjoint operator with discrete spectrum and
$\lam\in\RR$, we denote by $N(\lam,A)$ the number of the
eigenvalues of $A$ not exceeding $\lam$ (counting multiplicities).

Let $\bfB_a^\pm$ denote the restriction of $\bfB_a$ to the spaces
$\Gam(\tilW,\tilF^\pm)$. Similarly,  let $\mB_{\pm}, \mB_{\pm,a}$
denote the restriction of the operators $\mB, \mB_a$ to the spaces
$\Gam(M\times\RR,\tilF^\pm)$.
\prop{N}
Let $\lam_{\pm}$ denote the smallest non-zero eigenvalue of
$(\mB_{\pm})^2$. Then, for any $0<\eps<\min\{\lam_+,\lam_-\}$,
there exists $A=A(\eps,V)>0$, such that
\eq{N}
    N\big(\lam_{\pm}-\eps, (\bfB_a^{\pm})^2\big)
    \ = \ \dim\Ker(\mB_{\pm})^2,\qquad \text{for all}
    \quad a>A.
\end{equation}

\eprop
Before proving the proposition let us explain how it implies
\reft{cobinv}.

\subsection{Proof of \reft{cobinv}}\Label{SS:pra-a}
Let $V^{\pm}_{\eps,a}\subset \Gam(\tilW,\tilF^\pm)$ denote the
vector space spanned by the eigenvectors of the operator
$(\bfB^{\pm}_a)^2$ with eigenvalues smaller or equal to
$\lam_{\pm}-\eps$. The operator $\bfB^{\pm}_a$ sends
$V^{\pm}_{\eps,a}$ into $V^{\mp}_{\eps,a}$. It follows that
\[
    \dim\Ker\bfB_a^{+} \ - \
            \dim\Ker\bfB_a^{-}  \ = \
    \dim V^{+}_{\eps,a} \ - \
            \dim V^{-}_{\eps,a}.
\]
By \refp{N}, the right hand side of this equality equals
$\dim\Ker\mB_{+}-\dim\Ker\mB_{-}$. Thus $\ind\bfB_a= \ind\mB$.
\reft{cobinv} follows now from Lemmas~\ref{L:depofa} and
\ref{L:model}. \hfill$\square$

The rest of this section is occupied with the proof of \refp{N}.

\subsection{Estimate from above on $N(\lam_{\pm}-\eps,(\bfB_a^{\pm})^2)$}\Label{SS:above}
We will first show that
\eq{above}
       N(\lam_{\pm}-\eps,(\bfB_a^{\pm)^2}))\le\dim\Ker\mB_{\pm}.
\end{equation}
To this end we will estimate the operator $\H$ from below. We will
use the technique of   \cite{Sh5,BrFar1}, adding some necessary
modifications.

\subsection{The IMS localization}\Label{SS:IMS}
Let $j,\oj:\RR\to[0,1]$ be  smooth functions such
that $j^2+\oj^2\equiv0$ and $j(t)=1$ for $t\ge3$, while $j(t)=0$
for $t\le 2$. Set \/ $j_a(t)=j(a^{-1/2}t)$, \ \
$\oj_a(t)=\oj(a^{-1/2}t)$. This functions induce smooth functions
on the cylinder $M\times[0,1]$, which we denote by the same
letters. By a slight abuse of notation we will denote by the same
letters also the smooth functions on $\tilW$ given by the formulas
\/ \(
    j_a(x)=j(a^{-1/2}p(x)), \  \oj_a(x)=\oj(a^{-1/2}p(x)).
\)

The following version of IMS
\footnote{The abbreviation IMS stands for the
initials of R.~Ismagilov, J.~Morgan, I.~Sigal and B.~Simon}
localization formula is due to Shubin \cite[Lemma 3.1]{Sh5}.
\lem{sh1}The following operator identity holds
  \eq{sh1}
        \H=\oj_a \H\oj_a+j_a\H j_a+\frac12[\oj_a,[\oj_a,\H]]+\frac12[j_a,[j_a,\H]].
  \end{equation}
\elem
\prf
  Using the equality $j_a^2+\oj_a^2=1$ we can write
  $$
        \H=j_a^2\H+\oj_a^2\H=j_a\H j_a +\oj_a\H\oj_a+j_a[j_a,\H]+\oj_a[\oj_a,\H].
  $$
   Similarly, \/
  \(
        \H=\H j_a^2+\H\oj_a^2=j_a\H j_a +\oj_a\H\oj_a-[j_a,\H]j_a-[\oj_a,\H]\oj_a.
  \)
  \/ Summing these identities and dividing by 2, we come to \refe{sh1}.
\eprf
We will now estimate each of the summands in the right hand side
of \refe{sh1}.

\lem{sh3}
  There exists $A>0$, such that
  $\oj_a\H\oj_a\ge \frac{a^2}8\oj_a^2$, for all $a>A$.
\elem
\prf
   Note that $p(x)\le3 a^{1/2}$ for any $x$ in the support of
   $\oj_a$. Hence, if $a>36$, we have $\oj_a^2|p(x)-a|^2\ge
   \frac{a^2}4\oj_a^2$.

   Set
   \(
      A \ = \  \max\big\{\, 36, 4\sup_{x\in\tilW}|R|^{1/2}\, \big\}
   \)
   and let $a>A$. Using \refl{bfB2}, we obtain

   {\(\displaystyle\hskip1.8in
        \oj_a\H\oj_a \ \ge \oj_a^2 |p(x)-a|^2 \ - \ \oj_a R\oj_a
        \ \ge \ \frac{a^2}8\oj_a^2.
   \)}
\eprf
\subsection{}\Label{SS:Pa}
Let $P_{a}:L^2(M\times\RR,\tilF)\to \Ker\mB_a$ be the orthogonal
projection. Let $P^{\pm}_a$ denote the restriction of $P_a$ to the
space $L^2(M\times\RR,\tilF^\pm)$. Then $P^{\pm}_a$ is a finite
rank operator and its rank equals $\dim\Ker\mB_{\pm,a}$. Clearly,
\eq{ge}
        \mB_{\pm,a}+\lam_{\pm}P^{\pm}_a\ge \lam_{\pm}.
\end{equation}
By identifying the support of $j_a$ in $M\times\RR$ with a subset
of $\tilW$, we can and we will consider $j_aP_aj_a$ and
$j_a\mB_aj_a$ as operators on $\tilW$. Then $j_a\H
j_a=j_a\mB_aj_a$. Hence, \refe{ge} implies the following
\lem{local} \ \ \
  \(\displaystyle
        j_a\bfB_a^{\pm} j_a+\lam_{\pm}j_aP^{\pm}_aj_a\ge
        \lam_{\pm}j_a^2, \qquad
                \rk j_aP^{\pm}_aj_a\le \dim\Ker\mB_{\pm}.
  \)
\elem
For an operator $A:L^2(\tilW,\tilF)\to L^2(\tilW,\tilF)$, we
denote by $\|A\|$ its norm.
\lem{sh2}Let $\displaystyle
   C \ = \ 2\max\, \Big\{\, \max\{|j'(t)|^2,|\oj'(t)|^2\}: \,
        t\in\RR\, \Big\}$.
  Then
  \eq{sh2}
        \|[j_a,[j_a,\H]\|\le Ca^{-1}, \qquad   \|[\oj_a,[\oj_a,\H]\|\le
        Ca^{-1}, \qquad \text{for all} \quad a>0.
  \end{equation}
\elem
\prf
From \refl{bfB2} we obtain

 \(\quad
    |[j_a,[j_a,\H]|=2|j_a'(t)|^2=2a^{-1/2}|j'(a^{-1/2}t)|, \qquad
    |[\oj_a,[\oj_a,\H]|=2a^{-1/2}|j'(a^{-1/2}t)|.
 \)
\eprf

From Lemmas~\ref{L:sh1}, \ref{L:local} and \ref{L:sh2} we obtain
the following
\cor{Da>}
For any $\eps>0$, there exists $A=A(\eps,V)>0$, such that, for all
$a>A$, we have
\eq{Da>}
    \bfB_a^{\pm}+\lam_{\pm}j_aP_a^{\pm}j_a \ \ge \ \lam_{\pm}-\eps,
    \qquad
     \rk j_aP^{\pm}_aj_a\le \dim\Ker\mB_{\pm}.
\end{equation}
\ecor
The estimate \refe{above} follows  from \refc{Da>} and the
following general lemma \cite[p.~270]{ReSi4}:
\lem{general}
   Assume that $A, B$ are self-adjoint operators in a Hilbert space
   $\calH$ such that  $\rk B\le k$ and there exists $\mu>0$
   such that
   $
        \langle (A+B)u,u\rangle \ge \mu\langle u,u\rangle
                \quad \text{for any} \quad u\in\Dom(A).
   $
   Then $N(\mu -\eps, A)\le k$ for any $\eps>0$.
\elem

\subsection{Estimate from below on $N(\lam_{\pm}-\eps,(\bfB_a^{\pm})^2)$}\Label{SS:below}
To prove \refp{N} it remains now to show that
\eq{below}
       N(\lam_{\pm}-\eps,(\bfB_a^{\pm})^2) \ \ge \ \dim\Ker\mB_{\pm} \equiv
       \dim\Ker\mB_{\pm,a}.
\end{equation}
Let $V^{\pm}_{\eps,a}\subset L^2(\tilW,\tilF)$ denote the vector
space spanned by the eigenvectors of the operator
$(\bfB^{\pm}_a)^2$ with eigenvalues smaller or equal to
$\lam_{\pm}-\eps$. Let $\Pi^{\pm}_{\eps,a}:L^2(\tilW,\tilF^\pm)\to
V^{\pm}_{\eps,a}$ be the orthogonal projection.  Then \/
 \(
        \rk \Pi^{\pm}_{\eps,a} \ = \ N(\lam_{\pm}-\eps,(\bfB_a^{\pm})^2).
 \)
As in \refss{Pa}, we can and we will consider
$j_a\Pi^{\pm}_{\eps,a}j_a$ as an operator on
$L^2(M\times\RR,\tilF^\pm)$. The proof of the following lemma does
not differ from the proof of \refc{Da>}.
\lem{below}
For any $\eps>0$, there exists $A=A(\eps,V)>0$, such that, for any
$a>A$, we have
\eq{Da<}
    \mB_{\pm,a}+\lam_{\pm}j_a\Pi_a^{\pm}j_a \ \ge \ \lam_{\pm}-\eps,
    \qquad
     \rk j_a\Pi^{\pm}_aj_a\le \dim N(\lam_{\pm}-\eps,(\bfB_a^{\pm})^2).
\end{equation}
\elem
The estimate \refe{below} follows now from Lemmas~\ref{L:below}
and \ref{L:general}.

The proof of \refp{N} is complete. \hfil$\square$

\

\noindent\small{\textit{Acknowledgements.} \ Some of the ideas
used in this paper I have learned from John Roe. I would like to
thank him for very useful and stimulating discussion we had.}

\bibliographystyle{amsplain}
\providecommand{\bysame}{\leavevmode\hbox
to3em{\hrulefill}\thinspace}

\end{document}